\documentclass[12pt]{amsart}
\usepackage{amscd}

\newtheorem{thm}{Theorem}[section]

\theoremstyle{definition}

\def\ph{\phi}

\def \refeq#1{equation (\ref{#1})}
\def \ra{\rightarrow}

\def \hom{\mbox{\rm Hom}}

\def \tns{\otimes}

\def \k{\mbox{$\mathfrak K$}}

\def \Z{\mbox{$\mathbb Z$}}

\def \P{\mbox{$\mathbb P$}}
\def\br#1#2{\lbrack#1,#2\rbrack}

\def\zt{\mbox{$\Z_2$}}

\def\sh{\operatorname{Sh}}

\def\inv{^{-1}}
\def\d{d}
\def\td{\tilde\d}

\def\im{\operatorname{Im}}
\def\A{\mbox{$\mathcal A$}}
\def\B{\mbox{$\mathcal B$}}
\def\L{L}

\def\LA{\mbox{$\L_{\A}$}}
\def\LB{\mbox{$\L_{\B}$}}

\def\m{\mbox{$\mathfrak m$}}
\def\a{\mbox{$\mathfrak a$}}

\def\p{\epsilon}

\def\coder{\operatorname{Coder}}
\def\diag{\operatorname{diag}}

\def\linf{\mbox{$L_\infty$}}
\def\and{\mbox{ \rm and }}

\def\s#1{(-1)^{#1}}
\DeclareMathOperator*{\invlim}{\underleftarrow{\rm lim}}

\def\htns{\hat\tns}
\def\SW{S(W)}

\def\htns{\hat\tns}

\def\inv{^{-1}}

\def\dt#1{d_{#1}}
\def\dinfn{\dt{\infty,n}}
\def\dinfne{\dt{\infty,n,e}}
\def\thmref#1{Theorem (\ref{#1})}

\author{Alice Fialowski}
\address{E\"otv\"os Lor\'and University\\
Department of Applied Analysis\\
H-1117 Budapest, P\'azm\'any P. s\'et\'any. 1/C, Hungary}
\email{fialowsk@cs.elte.hu}
\author{Michael Penkava}
\address{University of Wisconsin\\
Department of Mathematics\\
Eau Claire, WI 54702-4004} \email{penkavmr@uwec.edu}
\subjclass{14D15,13D10,14B12,16S80,16E40,\\17B55,17B70}
\keywords{Strongly Homotopy Lie Algebras, $L_\infty$ Algebras, Superalgebras,
Lie Algebras, Extensions, Moduli space}
\thanks{The research of the authors was partially supported by
grants MTA-OTKA-NSF 38453, OTKA T043641, T043034, by the Humboldt
Foundation, and grants from the University of Wisconsin-Eau Claire. The
final version of the paper was completed during the stay of the first
author at the Max-Planck-Institut f\"ur Mathematik in Bonn.}
\title[$1|2$ dimensional \linf\ Algebras]%
{Strongly Homotopy Lie Algebras of One Even and Two Odd Dimension}

\begin{document}
\setlength{\multlinegap}{0pt}
%\nocite{ps2,pen1,ls,pen2,pen3,kon,
%fm,mar,mar2,ksv,sta1,sta2,sta3,getz,getz2,ge_ka1,ge_ka2,
%gers,lod,hoch,fi,fi2,ff2,ff,ff3}
\begin{abstract}
 We classify strongly homotopy Lie algebras - also called \linf\
 algebras -  of one even and two  odd dimension, which
are related to $2|1$-dimensional \zt-graded Lie algebras. What makes
this case interesting is that there are many nonequivalent \linf\
examples, in addition to the \zt-graded Lie algebra (or superalgebra)
structures, yet the moduli space is simple enough that we can give a
complete classification up to equivalence.
\end{abstract}
\date\today
\maketitle
%\table
%\input intro.tex
\vspace*{ -15pt}
\section{Introduction}
Strongly homotopy Lie algebras, or
  \linf\ algebras, for short, are natural
generalizations of Lie algebras and Lie superalgebras. In this paper,
all spaces will be \zt-graded, and we will work in the parity reversed
definition of the \linf\ structure. From this perspective,
an ordinary
$3$-dimensional Lie algebra is the same as an \linf\ algebra
structure on a $0|3$ (0 even and 3 odd) dimensional \zt-graded
vector space. Strongly
homotopy Lie algebras were first described in \cite{ss} and have
recently been the focus of much attention both in mathematics
(\cite{sta3,HSch}) and in mathematical physics
(\cite{aksz,bfls,ls, RW,M})

Even though \linf\ algebras appear in many contexts, there are
only a few examples known, even in low dimensions, other than
Lie superalgebras and differential graded Lie algebras.
In \cite{fp2},
the authors classified \linf\ algebras of dimension less than or
equal to 2. In \cite{fp3}, we constructed miniversal deformations
for all \linf\ structures on a space of dimension $0|3$. Since $0|3$-dimensional
vector spaces in our model of \linf\ algebras correspond
to $3|0$-dimensional spaces in the usual grading, the theory of
$0|3$-dimensional \linf\ algebras coincides with the theory of Lie algebra
structures on non-graded 3-dimensional vector spaces, whose
classification is well known. Because the symmetric algebra of a
$0|3$-dimensional vector space is finite dimensional, and there are no
odd cochains of any degree other than 2, the theory of \linf\ algebras
over a $0|3$-dimensional vector space introduces no new features beyond
the usual classification problem and deformation theory
for $3$-dimensional ordinary Lie algebras.

In \cite{bfp1}, \linf\ algebras of dimension $2|1$ are studied,
and a classification of the \linf\ algebra structures determined
by codifferentials whose leading term has degree less than or
equal to two is given. That classification includes the
deformation of \zt-graded Lie algebra structures into more general
\linf\ structures. The theory in the $2|1$ dimensional case is
more complicated than the $1|2$-dimensional algebras we will study
in this paper, because the space of $n$-cochains on a $1|2$
dimensional space has dimension $6|6$ for $n>1$, while the space
of $n$-cochains on a $2|1$- dimensional space has dimension
$3n+2|3n+1$, making it more difficult to classify the
nonequivalent structures.  Thus the $1|2$ dimensional case which
we study here is intermediate in complexity between the $0|3$ and the
$2|1$ dimensional cases.

Our purpose in this paper is to add to the body of examples of \linf\
algebras by classifying all such structures on a $1|2$-dimensional
vector space over an algebraically closed  field $\k$
of characteristic zero. We hope that these examples will help
to build a picture of the techniques necessary to attack the
classification problem in greater generality.

In Section~2 we introduce $L_\infty$ algebras, deformations, and give
some useful facts on equivalent codifferentials and extensions.
Section~3 treats the classification of codifferentials with leading
term $d_1\ne 0$. In Section~4 we classify codifferentials with leading
term $d_N$, where $N>1$. Sections 5 through 8 deal with cohomology
computations for the codifferentials we obtained in the previous
section, in order to get a better picture of the moduli space. In
Section 9 we determine extensions of each of the codifferentials.

\section{Basic Definitions and Facts}

\subsection{ Strongly Homotopy Lie Algebras}

A detailed introduction to \linf\ algebras can be obtained in
\cite{lm, ls, pen3, pen4}. Here we summarize the basic notions. If $W$ is
a \zt-graded vector space, then $\SW$ denotes the symmetric coalgebra
of $W$. If we let $T(W)$ be the reduced tensor algebra
$T(W)=\bigoplus_{n=1}^\infty W^{\tns n}$, then the reduced symmetric
algebra $S(W)$ is the quotient of the tensor algebra by the graded
ideal generated by $u\tns v-\s{uv}v\tns u$ for elements $u$, $v\in W$.
The symmetric algebra has a natural coalgebra structure, given by
\begin{equation*}
\Delta(w_1\dots w_n)=
\sum_{k=1}^{n-1}\sum_{\sigma\in\sh(k,n-k)}\epsilon(\sigma)
w_{\sigma(1)}\dots w_{\sigma(k)}\tns w_{\sigma(k+1)}\dots w_{\sigma(n)},
\end{equation*}
where we denote the product in $S(W)$ by juxtaposition, $\sh(k,n-k)$
is the set of \emph{unshuffles} of type $(k,n-k)$, and
$\epsilon(\sigma)$ is a sign determined by $\sigma$ (and $w_1\dots
w_n$) given by
$$
w_{\sigma(1)}\dots w_{\sigma(n)}=\epsilon(\sigma)w_1\dots w_n.
$$
A {\it coderivation} on $S(W)$ is a map $\delta:S(W)\ra S(W)$
satisfying $$\Delta\circ\delta=(\delta\tns
I+I\tns\delta)\circ\Delta.$$ Let us suppose that the even part of $W$
has basis $e_1\dots e_m$, and the odd part has basis $f_1\dots f_n$,
so that $W$ is an $m|n$ dimensional space. Then a basis of $\SW$ is
given by all vectors of the form $e_1^{k_1}\dots
e_m^{k_m}f_1^{l_1}\dots f_n^{l_n}$, where $k_i$ is any nonnegative
integer, and $l_i\in\zt$. An \linf\ structure on $W$ is simply an odd
codifferential on $\SW$, that is to say, an odd coderivation whose
square is zero. The space $\coder(W)$ can be naturally identified with
$\hom(S(W),W)$, and the Lie superalgebra structure on $\coder(W)$
determines a Lie bracket on $\hom(S(W),W)$ as follows. Denote
$L_m=\hom(S^m(W),W)$ so that $L=\hom(S(W),W)$ is the direct product of
the spaces $L_m$. If $\alpha\in L_m$ and $\beta\in L_n$, then
$\br\alpha\beta$ is the element in $L_{m+n-1}$ determined by
\begin{multline}\label{braform}
\br\alpha\beta
(w_1\dots w_{m+n-1})=\\
\sum_{\sigma\in\sh(n,m-1)}\epsilon(\sigma)\alpha(\beta(w_{\sigma(1)}\dots
w_{\sigma(n)})w_{\sigma(n+1)}\dots w_{\sigma(m+n-1)})\\
-\s{\alpha\beta}\
\sum_{\sigma\in\sh(m,n-1)}\epsilon(\sigma)
\beta(\alpha(w_{\sigma(1)}\dots
w_{\sigma(m)})w_{\sigma(m+1)}\dots w_{\sigma(m+n-1)}).
\end{multline}
Another way to express this bracket is in the form
\begin{equation*}
\br\alpha\beta=\alpha\tilde\beta-\s{\alpha\beta}\beta\tilde\alpha,
\end{equation*}
 where for $\ph\in\hom(S^k(W),W)$, $\tilde\ph$ is the
associated coderivation, given by
\begin{equation*}
\tilde\ph(w_1\dots w_n)=
\sum_{\sigma\in\sh(k,n-k)}\epsilon(\sigma)
\ph(w_{\sigma(1)}\dots w_{\sigma(k)})w_{\sigma(k+1)}\dots w_{\sigma(n)}.
\end{equation*}

If $W$ is completely odd, and $d\in\L_2$, then $d$ determines an ordinary
Lie algebra on $W$, or rather on its parity reversion.  This is the case
we considered in \cite{fp3}.  In that case, the
symmetric algebra on $W$ looks like the exterior algebra on $W$ if
we forget the grading.  If we define $[a,b]=d(ab)$ for $a,b\in W$, then
the bracket is antisymmetric because $ba=-ab$, and moreover
\begin{align*}
0=\br dd(abc)=&\
\frac12\sum_{\sigma\in\sh(2,1)}\epsilon(\sigma)
d(d(\sigma(a)\sigma(b))\sigma(c))\\=&
d((d(ab)c)+d(d(bc)a)-d(d(ac)b)\\
=&[[a,b],c]+[[b,c],a]-[[a,c],b],
\end{align*}
which is the Jacobi identity.  When $d\in\L_2$ and $W$ has a true
grading, then the same principle holds, except that one has to take
into account a sign arising from the map $S^2(W)\ra\bigwedge^2(V)$,
where $V$ is the parity reversion of $W$, so the induced bracket on
$V$ satisfies the graded Jacobi identity. Thus \zt-graded Lie
algebras are also examples of \linf\ algebras.

We work in the framework of the parity reversion $W=\Pi V$ of the
usual vector space $V$ on which an \linf\ algebra structure is
defined,  because in the $W$ framework  an \linf\ structure is simply
an odd coderivation $d$ of the symmetric coalgebra $S(W)$, satisfying
$d^2=0$,  in other words, it is an odd codifferential in the
\zt-graded Lie algebra of coderivations of $S(W)$.  As a consequence,
when studying \zt-graded Lie algebra structures on $V$, the parity is
reversed,  so that a $1|2$-dimensional vector space $W$ corresponds to
a $2|1$-dimensional \zt-graded Lie structure on $V$.  Moreover,  the
anti-symmetry of the Lie bracket on $V$ becomes \zt-graded symmetry of the
associated coderivation $d$ on $S(W)$.

A formal power series $d=d_1+\cdots$, with $d_i\in L_i$ determines an
element in $L=\hom(S(W),W)$, which is naturally identified with
$\coder(S(W))$, the space of coderivations of the symmetric coalgebra
$S(W)$. Thus $L$ is a \zt-graded Lie algebra, and we say that \emph{an
odd element $d$ in $L$ is a codifferential} if $d^2=0$, which is
equivalent to $\br dd=0$. Then $d$ is said to be an \linf\ structure
on $W$. When $d=d_2$, the \linf\ algebra is just a Lie superalgebra,
and when $d=d_1+d_2$, the $d_2$ cochain determines a Lie superalgebra
structure on $W$, while $d_1$ is a differential on the algebra, so that
$d$ is just a differential graded Lie algebra structure on $W$.

If $g=g_1+\cdots\in\hom(S(W),W)$, and $g_1:W\ra W$ is invertible,
then it determines a coalgebra automorphism of $S(W)$:
 $\tilde g:S(W)\ra S(W')$, that is a map satisfying
$$
\Delta'\circ\tilde g=(\tilde g\tns \tilde g)\circ\Delta.
$$
If $d$ and $d'$ are \linf\ algebra structures on $W$ and $W'$, resp.,
then $\tilde g$ is a homomorphism between these structures if $\tilde
g\circ d=\d'\circ\tilde g$. Two \linf\ structures $d$ and $d'$ on $W$
are equivalent, and we write $d'\sim d$, when there is a coalgebra
automorphism $\tilde g$ of $S(W)$ such that $d'={\tilde g}^*(d)={\tilde
g}\inv\circ d\circ {\tilde g}$.   Furthermore, if $d=d'$, then $\tilde
g$ is said to be an automorphism of the \linf\ algebra.

\subsection{Deformations}

An augmented local ring $\A$ with maximal ideal $\m$ will be called an
\emph{infinitesimal base} if  $\m^2=0$, and a  \emph{formal base} if
$\A=\invlim_n \A/\m^n$. A deformation of an \linf\ algebra structure
$d$ on $W$ with  base given by a local ring $\A$ with augmentation
$\epsilon:\A\ra\k$, where $\k$ is the field over which $W$ is defined,
is an $\A$-\linf\ structure $\td$ on $W\tns\A$ such that the morphism
of $\A$-\linf\ algebras $\epsilon_*=1\tns\epsilon:\LA=L\htns \A\ra
L\tns\k=L$ satisfies $\epsilon_*(\td)=d$.  (Here $L\htns\A$ is an
appropriate completion of $L\tns\A$.) The deformation is called
infinitesimal (formal) if $\A$ is an infinitesimal (formal) base.

In general, the cohomology $H(D)$ of $d$ given by the operator $D:L\ra
L$ with $D(\ph)=\br\ph d$ may not be finite dimensional, and does not
have a natural \Z-grading, as the cohomology of Lie algebras does. Moreover,
the entire odd part of the cohomology governs the deformations of the
\linf\ algebra structure, while the even part of the cohomology determines
infinitesimal automorphisms.
Nevertheless, $L$
has a natural filtration $L^n=\prod_{i=n}^\infty L_i$, which induces a
filtration $H^n$ on the cohomology, because $D$ respects the
filtration. We say that $H(D)$ is of finite type if $H^n/H^{n+1}$ is finite
dimensional. Since this is always true when $W$ is finite dimensional,
the examples we study here will be of finite type.  A set
$\delta_i$ will be called a \emph{basis of the cohomology}, if any
element $\delta$ of the cohomology can be expressed uniquely as a
formal sum $\delta=\delta_i a^i$. If we identify $H(D)$ with a
subspace of the space of cocycles $Z(D)$, and we choose a \emph{basis}
$\beta_i$ of the coboundary space $B(D)$, then any element
$\zeta\in\Z(D)$ can be expressed uniquely as a sum $\zeta=\delta_i a^i
+\beta_i b^i$.

For each $\delta_i$, let $u^i$ be a parameter of opposite parity.  Then
the infinitesimal deformation $d^1=d+\delta_i u^i$, with base
$\A=\k[u_i]/(u_iu_j)$ is universal in the sense that if $\td$ is any
infinitesimal deformation with base $\B$, then there is a unique
homomorphism $f:\A\ra\B$, such that the morphism $f_*=1\tns
f:\LA\ra\LB$ satisfies $f_*(d^1)\sim \td$.

For formal deformations, there is no universal object in the sense
above. A \emph{versal deformation} is a deformation $d^\infty$ with
formal base $\A$ such that if $\td$ is any formal deformation with base
$\B$, then there is some morphism $f:\A\ra\B$ such that
$f_*(d^\infty)\sim\td$. If $f$ is unique whenever $\B$ is
infinitesimal, then the versal deformation is called \emph{miniversal}.
In \cite{fp1}, we constructed a miniversal deformation for \linf\
algebras with finite type cohomology in general, and in \cite{fp2,fp3}
we computed versal deformations of \linf\ algebras on spaces of dimension
two or less.

In this paper we do not compute versal deformations, but we use the
concept of deformation in describing the moduli space of \linf\ structures
on our $1|2$-dimensional space $W$.

\subsection{Equivalent codifferentials and extensions}

 We will construct all equivalence classes of
\linf\ structures on $W$. We shall use the following facts, which
are established in \cite{fp1,bfp2}, to aid in the classification.

If $d$ is an \linf\ structure on $W$, and $d_N$ is the first
nonvanishing term in $d$, then $d_N$ is a codifferential, which we
call the \emph{leading term of $d$}, and we say that $d$ is
\emph{an extension} of $d_N$. Define the cohomology operator $D$
by $D(\ph)=\br\phi {d_N}$, for $\ph\in L$. Then the following
formula holds for any extension $d$ of $d_N$ as an \linf\
structure, and all $n\ge N$:
\begin{equation}\label{extension}
D(d_{n+1})=-\frac12\sum_{k=N+1}^{n}\br{d_k}{d_{n+N-k+1}}.
\end{equation}
Note that the terms on the right all have index less than $n+1$.
If a coderivation $d$ has been constructed up to terms of degree
$m$, satisfying \refeq{extension} for $n=1\dots m-1$, then the
right hand side of \refeq{extension} for $n=m$ is automatically a
cocycle. Thus $d$ can be extended to the next level precisely when
the right hand side is a coboundary. There may be many
nonequivalent extensions, because the term $d_{n+1}$ which we add
to extend the coderivation is only determined up to a cocycle.  An
extension $d$ is given by any coderivation which satisfies
\refeq{extension} for every $m=N+1\dots$ whose leading term is
$d_N$.

Classifying the extensions of $d_N$ can be quite complicated.
However, the following theorem often makes it easy to classify the
extensions.
\begin{thm}\label{zerocohomology}
If the cohomology $H^n(d_N)=0$, for $n>N$, then any extension of
$d_N$ to an \linf\ structure $d$ is equivalent to the structure
$d_N$.
\end{thm}

Before classifying the extensions of a codifferential $d_N$, we
need to classify the codifferentials in $L_N$ up to equivalence,
that is, we need to study the \emph{moduli space} of degree $N$
codifferentials. A \emph{linear automorphism} of $S(W)$ is an
automorphism determined by an isomorphism $g_1:W\ra W$. If $g$ is
an arbitrary automorphism, determined by maps $g_i:S^i(W)\ra W$,
and $W$ is finite dimensional, then $g_1$ is an isomorphism, so
this term alone induces an automorphism of $S(W)$ which we call
the \emph{linear part} of $g$.

The following theorem simplifies the classification of equivalence
classes of codifferentials in $L_N$.
\begin{thm}\label{simpleauto}
If $d$ and $d'$ are two codifferentials in $L_N$, and $g$ is an
equivalence between them, then the linear part of $g$ is also an
equivalence between them.
\end{thm}
Thus we can restrict ourself to linear automorphisms when determining
the equivalence classes of elements in $L_N$. If two codifferentials
are equivalent under a linear automorphism, then we say that they
are \emph{linearly equivalent}.
We will also use the following result.
\begin{thm}\label{simpleequiv}
Suppose that $d$ and $d'$ are equivalent codifferentials.  Then their
leading terms have the same degree and are linearly equivalent. \end{thm}

As a consequence of these theorems, we proceed to classify the
codifferentials as follows. First,  find all equivalence classes of
codifferentials of degree $N$. For each equivalence class, study the
equivalence classes of extensions of the codifferential.

Let us first establish some basic notation for the cochains. Suppose
$W=\langle w_1,w_2,w_3\rangle$, with $w_1$ and $w_2$ being odd
elements and $w_3$ an even element. An element in $W$ will be expressed
as a linear combination $w=w_1a+w_2b+w_3c$, with the coefficients written
on the right.
If $I=\{i_1,i_2,i_3\}$ is a
multi-index, with $i_1$ and $i_2$ either zero or one, let
$w_I=w_1^{i_1}w_2^{i_2}w_3^{i_3}$. For simplicity, we will denote
$w_I$ simply by $I$.

Then for $n\ge 0$,
\begin{equation*}S^{n+2}(W)=\langle (1,0,n+1),(0,1,n+1),(0,0,n+2),(1,1,n)
\rangle.
\end{equation*}

Note that $(1,0,n+1)$ and $(0,1,n+1)$ are odd, while $(0,0,n+2)$
and $(1,1,n)$ are even, so $\dim(S^{n+2}(W))=2|2$. If $g$ is a
linear automorphism of $S(W)$, then because automorphisms preserve
parity, we must have
\begin{align*}
g(w_1)=&w_1r+w_2p\\g(w_2)=&w_1r+w_2s\\g(w_3)=&w_3q,
\end{align*}
for some coefficients $l,p,r,s,q$ such that
$q(ls-pr)\ne 0$. If $q=1$, then we will express $g$ as a matrix
$g=\begin{pmatrix}l&r\\p&s\end{pmatrix}$. There is a unique
extension of $g$ to a linear automorphism of $S(W)$, given by
\begin{align*}
g(1,0,n+1)&=(1,0,n+1)lq^{n+1}+(0,1,n+1)pq^{n+1}\\
g(0,1,n+1)&=(1,0,n+1)rq^{n+1}+(0,1,n+1)sq^{n+1}\\
g(0,0,n+2)&=(0,0,n+2)q^{n+2}\\
g(1,1,n)&=(1,1,n)(ls-pr)q^n
\end{align*}

Let $L_n:=\hom(S^n(W),W)$. Define
\begin{equation*}
\ph^I_j(w_J)=I!\delta^I_J w_j,
\end{equation*}
where $I!=i_1!i_2!i_3!$. If we let  $|I|=i_1+i_2+i_3$, then
$L_n=\langle \ph^I_j, |I|=n \rangle$. If $\ph$ is odd, we denote
it by the symbol $\psi$ to make it easier to distinguish the even
and odd elements. Then
\begin{align*}
(L_{n+2})_0&=\langle\ph^{1,0,n+1}_1,\ph^{1,0,n+1}_2,\ph^{0,1,n+1}_1,
\ph^{0,1,n+1}_2,\ph^{0,0,n+2}_3,\ph^{1,1,n}_3\rangle\\
(L_{n+2})_1&=
\langle\psi^{1,0,n+1}_3,\psi^{0,1,n+1}_3,\psi^{0,0,n+2}_1,\psi^{0,0,n+2}_2,
\psi^{1,1,n}_1,\psi^{1,1,n}_2\rangle,
\end{align*}
so $\dim L_{n+2}=6|6$.

For $L_1$, $n=-1$, and no terms of the form $\{1,1,n\}$ occur, so
$\dim L_1=5|4$. We study the special case of codifferentials with
leading term in $L_1$ first.

\section{Classification of codifferentials with $d_1\ne0$}

Let us suppose that $d$ is an odd element of $L_1$. Then
\begin{equation*}
d=\psi^{1,0,0}_3a_1+\psi^{0,1,0}_3a_2+\psi^{0,0,1}_1a_3+\psi^{0,0,1}_2a_4.
\end{equation*}
We have
\begin{align*}
\frac12\br dd&= \ph^{1,0,0}_1a_1a_3 +\ph^{0,1,0}_1a_3a_2
+\ph^{1,0,0}_2a_4a_1 +\ph^{0,1,0}_2a_2a_4
\\&+\ph^{0,0,1}_3(a_2a_4+a_1a_3).
\end{align*}
For this bracket to vanish, it is clear that either $a_1=a_2=0$ or
$a_3=a_4=0$. Let us suppose that the former occurs, so that
$d=\psi^{0,0,1}_1a_3+\psi^{0,0,1}_2a_4$. Then we claim that $d$ is
equivalent to $d'=\psi^{0,0,1}_1$. For if we let $g$ be a linear
automorphism such that $g(w_1)=w_1a_3 +w_2 a_4$, then $g\inv d g=d'$.
Similarly, if $a_3=a_4=0$, then
$d=\psi^{1,0,0}_3a_1+\psi^{0,1,0}_3a_2$, which is equivalent to
$\ph^{0,1,0}_3$.  Thus there are exactly two equivalence classes of
codifferentials of degree 1.

In order to classify their extensions, we need to compute the
cohomology. Let $z_n=\dim(Z^n=\ker D:L_n\ra L_n)$ be the dimension of
the cocycles, and $b_n=\dim(B^n=\im D:L_n\ra L_n)$ be the dimension of
the space of coboundaries. Then $h_n=z_n-b_n$. (Note that dimensions
are represented as a pair $k|l$ consisting of the dimension of the
even, followed by the dimension of the odd part of the space.) Note
that for a codifferential of degree 1, $D:L_n\ra L_n$. It follows that
$h_1=\dim(H^1)$ cannot be zero. We shall see that for both
codifferentials, $h_n=\dim(H^n)=0$ for $n>1$.  Thus we can apply
\thmref{zerocohomology} to conclude that all codifferentials which
have a nonzero leading term of degree 1 are equivalent to one of these
two codifferentials.

Let $d_0=\psi^{1,0,0}_3$. The following table of coboundaries
allows us to compute $H^1$ immediately.
\begin{align*}
D(\ph^{1,0,0}_1)&=-\psi^{1,0,0}_3, &D(\psi^{1,0,0}_3)&=0\\
D(\ph^{1,0,0}_2)&=0, &D(\psi^{0,1,0}_3)&=0\\
D(\ph^{0,1,0}_1)&=-\psi^{0,1,0}_3, &D(\psi^{0,0,1}_1)&=\ph^{1,0,0}_1+\ph^{0,0,1}_3\\
D(\ph^{0,1,0}_2)&=0,&D(\psi^{0,0,1}_2)
    &=\ph^{1,0,0}_2\\
D(\ph^{0,0,1}_3)&=\psi^{1,0,0}_3\\
\end{align*}
Thus $z_1=3|2$,  $b_1=2|2$ and so $h_1=1|0$.  Moreover, we can take
$H^1=\langle \ph^{0,1,0}_2\rangle$. Since $\alpha=\ph^{0,1,0}_2a$ is even,
we can interpret this cohomology as a reflection of the fact that $\tilde g=\exp(\alpha)$,
which is a linear automorphism given by $g_1=\diag(1,e^a,1)$, leaves $d$ invariant. Actually,
$\exp(\alpha)$ leaves $d$ invariant for any even cocycle $\alpha$.  The automorphisms determined
by coboundaries are called \emph{inner automorphisms}, so we can interpret the fact that $h_1$
has even dimension 1 as the fact that the group of \emph{outer automorphisms} of
the codifferential $d$ has dimension 1.

Next, we compute coboundaries in
$L_{n+2}$ for $n\ge 0$.

\begin{align*}
D(\ph^{1,0,n+1}_1)&=-\psi^{1,0,n+1}_3, &D(\psi^{1,0,n+1}_3)&=0\\
D(\ph^{1,0,n+1}_2)&=0, &D(\psi^{0,1,n+1}_3)&=\ph^{1,1,n}_3(n+1)\\
D(\ph^{0,1,n+1}_1)&=\psi^{1,1,n}_1(n+1)-\psi^{1,0,n+1}_3,&
D(\psi^{0,0,n+2}_1)&=\ph^{1,0,n+1}_1(n+2)\\&&&+\psi^{0,0,n+2}_3\\
D(\ph^{0,1,n+1}_2)&=\psi^{1,1,n}_2(n+1), &
D(\psi^{0,0,n+2}_2)&=\ph^{1,0,n+1}_2(n+2)\\
D(\ph^{0,0,n+2}_3)&=\psi^{1,0,n+1}_3(n+2),&D(\psi^{1,1,n}_1)&=\ph^{1,1,n}_3\\
D(\ph^{1,1,n}_3)&=0,&D(\psi^{1,1,n}_2)&=0\\
\end{align*}
Notice that if you ignore the three input terms with negative indices, this
table coincides with the previous table when $n=-1$. From the table,
we can see that the space of even cocycles has dimension 3, so the
space of odd coboundaries must also have dimension 3. Similarly, the
odd cocycles span a three dimensional space, so the even coboundaries
span a three dimensional space. Thus $z_n=3|3=b_n$, so $h_n=0$ when
$n>1$.

%Let us compute the miniversal deformation of $d_0$.  We have
%\begin{equation*}
%d^1_0=\psi^{0,1,0}_3+\ph^{1,0,0}_2\theta,
%\end{equation*}
%where $\theta$ is an odd parameter.  Since $\theta^2=0$ automatically,
%we observe that $\br{d^1}{d^1}=0$, so the miniversal deformation
%coincides with the universal infinitesimal deformation, and there are
%no relations on the base.

If $d_*=\psi^{0,0,1}_1$, then a similar calculation gives
$H^1=\langle \ph^{0,1,0}_2\rangle$, and $H^n=0$ if $n>1$.
%Thus the universal infinitesimal and miniversal deformation is
    %$d^1_*=\ph^{0,0,1}_1+\ph^{0,1,0}_2\theta$.
So far, nothing much interesting has emerged.

\section{Codifferentials with leading term $d_N$, where $N>1$}

Suppose that $d$ is a codifferential in $L_{m+2}$, where $m\ge 0$.
Then we can express
\begin{equation*}
d=\psi^{1,0,m+1}_3a_1+\psi^{0,1,m+1}_3a_2+\psi^{0,0,m+2}_1a_3
    +\psi^{0,0,m+2}_2a_4+\psi^{1,1,m}_1a_5+\psi^{1,1,m}_2a_6.
\end{equation*}
The bracket of $d$ with itself is given by
\begin{align*}
\tfrac12\br dd&= \ph^{1,0,2m+2}_1(-a_4a_5+a_1a_3(m+2))
+\ph^{0,1,2m+2}_1a_3(a_5+a_2(m+2))\\
&+\ph^{1,0,2m+2}_2a_4(-a_6+a_1(m+2))
+\ph^{0,1,2m+2}_2(a_3a_6+a_2a_4(m+2))\\
&+\ph^{1,1,2m+1}_3(a_2a_6+a_1a_5)
+\ph^{0,0,2m+3}_3(a_2a_4+a_1a_3)\\
\end{align*}
A little arithmetic reduces the set of equations which follow from
the above to
\begin{align*}
a_2a_6+a_1a_5&=0\\
a_1a_3+a_2a_4&=0\\
a_3(a_5+(m+2)a_2)&=0\\
a_4(a_5+(m+2)a_2)&=0\\
a_3(a_6-(m+2)a_1)&=0\\
a_4(a_6-(m+2)a_1)&=0.
\end{align*}
It follows that either $a_3=a_4=0$ or both $a_5=-(m+2)a_2$ and
$a_6=(m+2)a_1$. Let us consider the case $a_3=a_4=0$ first. Then
the equations reduce to $a_2a_6+a_1a_5=0$. If $a_1=a_2=0$ then we
obtain no conditions on $a_5$ or $a_6$, so we obtain the
codifferential $d=\psi^{1,1,m}_1a_5+\psi^{1,1,m}_2a_6$. Let $b_1$,
$b_2$ be such that $b_2a_5-b_1a_6=1$. Let $d'=\psi^{1,1,m}_1$.
Then $d'=g\inv dg$ for the automorphism
$g=\begin{pmatrix}a_5&b_1\\a_6&b_2
\end{pmatrix}$, so $d'$ and $d$ are equivalent.

If either $a_1\ne 0$ or $a_2\ne 0$, then we have $a_5=ka_2$ and
$a_6=-ka_1$ for some $k$.  Thus we obtain the codifferential
\begin{equation*}
d=\psi^{1,0,m+1}_3a_1+\psi^{0,1,m+1}_3a_2+\psi^{1,1,m}_1ka_2
    -\psi^{1,1,m}_2ka_1.
\end{equation*}
Let $g=\begin{pmatrix}a_2&b_1\\-a_1&b_2\end{pmatrix}$, where
$a_1b_1 +a_2b_2=1$. Then one can easily check that $gd'=dg$, for
$d'=\psi^{0,1,m+1}_3+\psi^{11m}_1k$.

Finally,  suppose that at least one of $a_3$ or $a_4$ does not
vanish.  Then $a_5=-(m+2)a_2$ and $a_6=(m+2)a_1$. Define
$g=\begin{pmatrix}a_3&b_1\\a_4&b_2\end{pmatrix}$, where
$a_3b_2-a_4b_1=1$.  Let $k=(a_1b_1+a_2b_2)$. Then $d$ is
equivalent to the codifferential
$d'=\psi^{0,1,m+1}_3k+\psi^{0,0,m+2}_1-\psi^{1,1,m}_1k(m+2)$. If
$k\ne0$ then we can  reduce this codifferential to
$d=\psi^{0,1,m+1}_3+\psi^{0,0,m+2}_1-\psi^{1,1,m}_1(m+2)$.

\smallskip

Thus we arrive at the following table of codifferentials of degree
$m+2$:
\begin{align*}
d_\infty&=\psi^{1,1,m}_1\\
\dt\lambda&=\psi^{0,1,m+1}_3+\psi^{1,1,m}_1\lambda\\
d_*&=\psi^{0,0,m+2}_1\\
d_{\#}&=\psi^{0,1,m+1}_3+\psi^{0,0,m+2}_1-\psi^{1,1,m}_1(m+2)
\end{align*}
If we wish to emphasize the degree of the codifferential,  we will
indicate it parenthetically.  Thus $\dt\lambda(k)$ represents the
degree $k$ version of $\dt\lambda$. Note that for $m=-1$, only
$\dt0$ and $d_*$ are well defined, and they coincide with the
codifferentials of degree 1 which we have already studied. It is
easy to check that none of these four types of codifferentials is
equivalent to any of the others. Moreover, $\dt\lambda(k)$ is
equivalent to $\dt\mu(k)$ precisely when $\lambda=\mu$. Thus we
have classified all codifferentials in $L_{m+2}$ for $m\ge0$.

Notice that any neighborhood of $d_*$ contains an element
equivalent to $d_{\#}$ and similarly, any neighborhood of
$d_\lambda$, for $\lambda=-(m+2)$ also contains an element
equivalent to $d_{\#}$. We say that $d_{\#}$ is infinitesimally
close to these elements, or, in the terminology of M. Gerstenhaber,
 that they are jump deformations
of  $d_{\#}$  It seems reasonable to characterize the
moduli space of codifferentials of degree $N>1$ as consisting of
$\k\P^1$ plus two special points.  In order to get a better
picture of our moduli space,  we need to look at the cohomology.

\section{Cohomology of $d_\lambda=\psi^{0,1,m+1}_3+\psi^{1,1,m}_1\lambda$}

First, we compute a table of coboundaries for
$D=[\cdot,d_\lambda]$. The coboundaries of the even elements are
given by
\begin{align*}
&D(\ph^{1,0,n+1}_1)=-\psi^{1,1,m+n+1}_1(n+1)\\
&D(\ph^{1,0,n+1}_2)=\psi^{1,1,m+n+1}_2(\lambda-n-1)
-\psi^{1,0,m+n+2}_3\\
&D(\ph^{0,1,n+1}_1)=0\\
&D(\ph^{0,1,n+1}_2)=-\psi^{1,1,m+n+1}_1\lambda-\psi^{0,1,m+n+2}_3\\
&D(\ph^{0,0,n+2}_3)=-\psi^{1,1,m+n+1}_1\lambda m+\psi^{0,1,m+n+2}_3(n+1-m)\\
&D(\ph^{1,1,n}_3)=0\\
\end{align*}
The odd coboundaries are
\begin{align*}
&D(\psi^{1,0,n+1}_3)=-\ph^{1,1,m+n+1}_3(-\lambda-m+n)\\
&D(\psi^{0,1,n+1}_3)=0\\
&D(\psi^{0,0,n+2}_1)=\ph^{0,1,m+n+2}_1(\lambda+n+2)\\
&D(\psi^{0,0,n+2}_2)=-\ph^{1,0,m+n+2}_1\lambda+\ph^{0,1,m+n+2}_2(n+2)
    +\ph^{0,0,m+n+3}_3\\
&D(\psi^{1,1,n}_1)=0\\
&D(\psi^{1,1,n}_2)=\ph^{1,1,m+n+1}_3
\end{align*}

Let us define
\begin{align*}
&z_n=\ker(D:L_n\ra L_{n+m+1})\\
&b_n=\im(D:L_n\ra L_{n+m+1})
\end{align*}
Then $h_{n+2}=z_{n+2}-b_{n-m+1}$. Note that if $n<m$,
$h_{n+2}=z_{n+2}$,  while $h_{m+2}=z_{m+2}-b_1$. Furthermore,
$z_{n+2}+b_{n+2}=6|6$, for $n\ge 0$ and  $z_1+b_1=5|4$.

Let us first consider the case $n\ne -1$. Then there are exactly 3 even
independent $(n+2)$-cocycles, and therefore three odd independent
$(n+2)$-coboundaries. The situation for odd cocycles is a bit more
complicated. It is clear that generically, we have 3 odd
$(n+2)$-cocycles and therefore 3 even $(n+2)$-coboundaries. When
$\lambda=-(n+2)$, we obtain an extra cocycle $\psi^{0,0,n+2}_1$. Thus,
when $\lambda$ is a negative integer, the cohomology is a bit different
than in the generic case.

When $n=-1$, there are always three even independent cocycles, and
thus 2 odd independent coboundaries. Generically, there is only
one odd 1-cocycle, but when $\lambda=-(m+1)$, we pick up an extra
odd cocycle $\psi^{1,0,0}_3$. Notice, that in this case,
$\lambda$ is also a negative integer.

\subsection{Cohomology when $\lambda$  is not a negative integer}

Let us suppose that $\lambda$ is generic, that is, it is not a
negative integer. Then  $z_{n+2}=3|3$ and $b_{n+2}=3|3$ for $n\ge
0$. Also we have $z_1=3|1$, so $b_1=3|2$.  Therefore we have
\begin{align*}
&h_1=3|1\\
&h_{n+2}=3|3,\quad\text{ if $0\le n<m$}\\
&h_{m+2}=0|1\\
&h_{n+2}=0,\quad\text{if $n>m$}
\end{align*}
We can give a basis for the cohomology in each dimension.
\begin{align*}
&H^1=\langle \psi^{0,1,0}_3, \ph^{1,0,0}_1,\ph^{0,1,0}_1,
m\ph^{0,1,0}_2-\ph^{0,0,1}_3\rangle\\
&(H^{n+2})_o=\langle \psi^{0,1,n+1}_3,\psi^{1,1,n}_1,
\psi^{1,0,n+1}_3+\psi^{1,1,n}_2(-\lambda-m+n)\rangle\\
&(H^{n+2})_e=\langle\ph^{0,1,n+1}_1,\ph^{1,1,n}_3,
-\ph^{1,0,n+1}_1\lambda+\ph^{0,1,n+1}_2(n+1-m)+\ph^{0,0,n+2}_3
\rangle\\
&\quad\text{if $0\le n<m$}\\
&H^{m+2}=\langle \psi^{1,1,m}_1\rangle\\
&H^{n+2}=0\quad\text{if $n>m$}.
\end{align*}

The terms of degree $n+2$, with $-1\le n<m$, govern deformations
of $\dt\lambda$ into codifferentials of lower degree. Looking at
$H^1$, we see that $d_\lambda$ generically deforms into a
codifferential which is equivalent to the codifferential
$d_0(1)=\psi^{0,1,0}_3$ of degree 1.

When the leading term of a deformation is of degree $n+2$, where
$0\le n<m$, our codifferential deforms into three distinct types
of codifferentials of degree $n+2$.  First, we have
$\dt0=\psi^{0,1,n+1}_3$, then $\dt\infty=\psi^{1,1,n}_1$, and
finally
$\dt{-\lambda-m+n}=\psi^{1,0,n+1}_3-\psi^{1,1,n}_2(-\lambda-m+n)$.
Note that since $\lambda$ is not a negative integer,
$-\lambda-m+n\ne 0$, so we do not obtain two different ways to
deform into $\dt0$. Of course, we could deform into a sum of
terms of these types as well.

%We will
%see which combinations are possible when we compute the versal
%deformations of $\dt\lambda$ later.

We are especially interested in the deformations of $\dt\lambda$
on the same level, because they tell us how the moduli space of
degree $m+2$ codifferentials fits together.  Consider the degree
$m+2$ part of the universal infinitesimal deformation of
$\dt\lambda$.  It is of the form
\begin{equation*}
d_\lambda^1=\psi^{0,1,m+1}_3+\psi^{1,1,m}_1(\lambda+t).
\end{equation*}
Thus we see that $\dt\lambda$ deforms into $\dt{\lambda+t}$.  This
is a deformation along the family $\dt\lambda$.

Finally, note that because $h_{n+2}=0$ for $n>m$, there are no
nonequivalent extensions of $d_\lambda$.

\subsection{Cohomology for $\lambda=-p$, if  $p\ne m+1$}

Suppose that $\lambda=-p$ is a negative integer, not equal to
$-(m+1)$. Then, depending on the value of $p$, we have several
distinct cases.

\subsubsection{Case $p=1$}

Now $h_1=3|2$, and we have an extra even $(m+2)$-cohomology class,
$\ph^{1,1,m}_3$, so $h_{m+2}=1|1$. We obtain a deformation in the
$d_*(1)=\psi^{0,0,1}_1$ direction. Thus the table of cohomology
for $\dt\lambda$ is modified by
\begin{align*}
&H^1=\langle \psi^{0,1,0}_3,
\psi^{0,0,1}_1,\ph^{1,0,0}_1,\ph^{0,1,0}_1,
m\ph^{0,1,0}_2-\ph^{0,0,1}_3\rangle\\
&H^{m+2}=\langle \psi^{1,1,m}_1,\ph^{1,1,m}_3\rangle
\end{align*}

\subsubsection{Case $2\le p <m+2$}

In this case we obtain an extra odd $p$-cohomology class
$\psi^{0,0,p}_1$, which means we again obtain an extra deformation
in the $d_*(p)$ direction. Of course, there is a corresponding
extra even $(m+p+1)$-cohomology class $\ph^{0,1,m+p}_1$. Notice
that this value of $\lambda$ coincides with the value of
$\lambda$ at which we already noticed something peculiar happening
in the moduli space. The table of cohomology for $\dt\lambda$
becomes modified by
\begin{align*}
&(H^{p})_o=\langle \psi^{0,1,p-1}_3,\psi^{1,1,p-2}_1,
\psi^{1,0,p-1}_3-\psi^{1,1,p-2}_2(2p-m),\psi^{0,0,p}_1\rangle\\
&H^{m+p+1}=\langle\ph^{0,1,m+p}_1\rangle.
\end{align*}

\subsubsection{Case $p=m+2$}

In this case we get a deformation of $d_\lambda$  in the
$d_*=\psi^{0,0,m+2}_1$ direction on the same level, and an extra
even $(2m+3)$-cohomology class $\ph^{0,1,2m+2}_1$.  Thus we have
an extra direction of deformation in the moduli space of degree
$m+2$ codifferentials. The table of cohomology is modified by
\begin{align*}
&H^{m+2}=\langle \psi^{1,1,m}_1,\psi^{0,0,m+2}_1\rangle\\
&H^{2m+3}=\langle\ph^{0,1,2m+2}_1\rangle.
\end{align*}

\subsubsection{Case $p>m+2$}

In this final case, we obtain something really interesting,
because there is a nontrivial extension $\dt{\lambda,e}=d_\lambda
+ \psi^{0,0,n+2}_1$ of $\d_\lambda$. We will discuss this
extension later. The table of cohomology is altered in degrees $p$
and $m+p+1$ as follows:
\begin{align*}
&H^{p}=\langle \psi^{0,0,p}_1\rangle\\
&H^{m+p+1}=\langle\ph^{0,1,m+p}_1\rangle.
\end{align*}

\subsection{Cohomology for $\lambda=-(m+1)$}

For $n=-1$, there is an extra odd 1-cohomology class,
$\psi^{1,0,0}_3$, so there is an extra even $(m+2)$-cohomology
class $\ph^{1,1,m}_3$. Then $h_1=3|2$ and $h_{m+2}=1|1$. We also
obtain an extra odd $(m+1)$-cohomology class $\psi^{0,0,m+1}_1$
and an extra even $(2m+2)$-cohomology class $\ph^{0,1,2m+1}_1$.
Thus the table of cohomology becomes modified as follows:
\begin{align*}
&H^1=\langle \psi^{0,1,0}_3,\psi^{1,0,0}_3,
\ph^{1,0,0}_1,\ph^{0,1,0}_1,
m\ph^{0,1,0}_2-\ph^{0,0,1}_3\rangle\\
&(H^{m+1})_o=\langle \psi^{0,1,m}_3,\psi^{1,1,m-1}_1,
\psi^{1,0,m}_3-\psi^{1,1,m-1}_2(2m)\rangle\\
&H^{m+2}=\langle \psi^{1,1,m}_1,\ph^{1,1,m}_3\rangle\\
&H^{2m+2}=\langle\ph^{0,1,2m+2}_1\rangle.
\end{align*}

\section{Cohomology of $d_\infty=\psi^{1,1,m}_1$}

The coboundaries are
\begin{align*}
D(\ph^{1,0,n+1}_1)&=0,\qquad&D(\psi^{1,0,n+1}_3)&=\ph^{1,1,m+n+1}_3\\
D(\ph^{1,0,n+1}_2)&=\psi^{1,1,m+n+1}_2,\qquad&D(\psi^{0,1,n+1}_3)&=0\\
D(\ph^{0,1,n+1}_1)&=0,\qquad&D(\psi^{0,0,n+2}_1)&=\ph^{0,1,m+n+2}_1\\
D(\ph^{0,1,n+1}_2)&=-\psi^{1,1,m+n+1}_1,\qquad&D(\psi^{0,0,n+2}_2)
    &=-\ph^{1,0,m+n+2}_1\\
D(\ph^{0,0,n+2}_3)&=-\psi^{1,1,m+n+1}_1m,\qquad&D(\psi^{1,1,n}_1)&=0\\
D(\ph^{1,1,n}_3)&=0,\qquad&D(\psi^{1,1,n}_2)&=0\\
\end{align*}
It is easy to see that $z_1=3|1$ and $z_{n+2}=4|3$ for $n>-1$.
Thus $b_{n+2}=3|2$ for $n\ge -1$. From this we get
\begin{align*}
&h_1=3|1\\
&h_{n+2}=4|3,\quad\text{ if $0\le n<m$}\\
&h_{n+2}=1|1,\quad\text{if $n\ge m$}.
\end{align*}
Consequently we always have nontrivial extensions of $\dt\infty$. If we
look at $H^1$, we see that the odd cohomology class
$d_0=\psi^{0,1,0}_3$ will deform into a codifferential of type $d_0$ in
degree 1. For $0\le n<m$, we have three independent odd cohomology
classes, $d_0=\psi^{0,1,n+1}_3$, $\td\infty=\psi^{1,1,n}_1$ and
$\psi^{1,1,n}_2$. The last cohomology class also is of type $d_\infty$.
For $n\ge m$, there is only one odd cohomology class,
$\psi^{0,1,n+1}_3$. When $n=m$, this cohomology class governs
deformations along the family. For $n>m$,  the existence of this
cohomology class implies that there is a nontrivial extension of
$\d_\infty$. In fact, there are many nonequivalent extensions, which we
will discuss later.  A table of the basis of the cohomology for
$\dt\infty$ is
\begin{align*}
&H^1=\langle \psi^{0,1,0}_3, \ph^{1,0,0}_1,\ph^{0,1,0}_1,
\ph^{0,1,0}_2m-\ph^{0,0,1}_3\rangle\\
&(H^{n+2})_o=\langle \psi^{0,1,n+1}_3,\psi^{1,1,n}_1,
\psi^{1,1,n}_2\rangle\\
&(H^{n+2})_e=\langle\ph^{0,1,n+1}_1,\ph^{1,1,n}_3,
\ph^{1,0,n+1}_1,\ph^{0,1,n+1}_2m-\ph^{0,0,n+2}_3
\rangle\\
&\quad\text{if $0\le n<m$}\\
&H^{n+2}=\langle
\psi^{0,1,n+1}_3,\ph^{1,0,n+1}_2m-\ph^{0,0,n+2}_3\rangle
\quad\text{if $n\ge m$}.
\end{align*}

\section{Cohomology of $d_*=\psi^{0,0,m+2}_1$}

First, we compute a table of coboundaries for
$D=[\cdot,d_\lambda]$. The boundaries of the elements are given by

\begin{align*}
&D(\ph^{1,0,n+1}_1)=\psi^{0,0,m+n+3}_1,\quad
    D(\psi^{1,0,n+1}_3) =\ph^{1,0,m+n+2}_1(m+2)+\ph^{0,0,m+n+3}_3\\
&D(\ph^{1,0,n+1}_2)=\psi^{0,0,m+n+3}_2,\qquad\qquad
    D(\psi^{0,1,n+1}_3) =\ph^{0,1,m+n+2}_1(m+2)\\
&D(\ph^{0,1,n+1}_1)=0,\qquad\qquad\qquad\qquad\qquad\qquad\quad
    D(\psi^{0,0,n+2}_1)=0\\
&D(\ph^{0,1,n+1}_2)=0,\qquad\qquad\qquad\qquad\qquad\qquad\quad
    D(\psi^{0,0,n+2}_2)=0\\
&D(\ph^{0,0,n+2}_3)=-\psi^{0,0,m+n+3}_1(m+2),\qquad\qquad\quad\
    -D(\psi^{1,1,n}_1)=\ph^{0,1,m+n+2}_1\\
&D(\ph^{1,1,n}_3)=-\psi^{1,1,m+n+1}_1(m+2)+\psi^{0,1,m+n+2}_3,\quad
    D(\psi^{1,1,n}_2)=\ph^{0,1,m+n+2}_2\\
\end{align*}
It is easy to compute that
\begin{align*}
&h_1=3|2\\
&h_{n+2}=3|3,\quad\text{ if $0\le n<m$}\\
&h_{m+2}=1|1\\
&h_{n+2}=0,\quad\text{if $n> m$}.
\end{align*}
A table of the cohomology is given by
\begin{align*}
&H^1=\langle \psi^{0,0,1}_1,\psi^{0,0,1}_2,
\ph^{0,1,0}_1,\ph^{0,1,0}_2,
\ph^{1,0,0}_2(m+2)+\ph^{0,0,1}_3\rangle\\
&(H^{n+2})_o=\langle \psi^{0,0,n+2}_1,\psi^{0,0,n+2}_2,
\psi^{0,1,n+1}_3-\psi^{1,1,n}_1(m+2)\rangle\\
&(H^{n+2})_e=\langle\ph^{0,1,n+1}_1,\ph^{0,1,n+1}_2,
\ph^{1,0,n+1}_2(m+2)+\ph^{0,0,n+2}_3
\rangle\\
&\quad\text{if $0\le n<m$}\\
&H^{m+2}=\langle\psi^{0,1,m+1}_3-\psi^{1,1,m}_1(m+2),\ph^{0,1,m+1}_2\rangle\\
&H^{n+2}=0,\quad \text{if $n>m$}.
\end{align*}
In degree 1, our codifferential deforms into $d_*$ in two
different ways. For $0\le n<m$, we have two deformations in the
$d_*$ direction, but we also get a deformation into the special
element $d_{-(m+2)}$ of the family.  This direction of deformation
is the only one that remains in degree $m+2$, and there are no
higher order deformations, thus no nontrivial extensions of $d_*$.

\section{Cohomology of $d_{\#}=\psi^{0,1,m+1}_3+\psi^{0,0,m+2}_1
    -\psi^{1,1,m}_1(m+2)$}

The even coboundaries of the elements are
\begin{align*}
D(\ph^{1,0,n+1}_1)&=-\psi^{1,1,m+n+1}_1(n+1)+\psi^{0,0,m+n+3}_1\\
D(\ph^{1,0,n+1}_2)&=-\psi^{1,1,m+n+1}_2(n+m+3)+\psi^{0,0,m+n+3}_2
    -\psi^{1,0,m+n+2}_3\\
D(\ph^{0,1,n+1}_1)&=0\\
D(\ph^{0,1,n+1}_2)&=\psi^{1,1,m+n+1}_1(m+2)-\psi^{0,1,m+n+2}_3\\
D(\ph^{0,0,n+2}_3)&=-\psi^{1,1,m+n+1}_1m(m+2)-\psi^{0,0,m+n+3}_1(m+2)\\
    &+\psi^{0,1,m+n+2}_3(n+1-m)\\
D(\ph^{1,1,n}_3)&=-\psi^{1,1,m+n+1}_1(m+2)+\psi^{0,1,m+n+2}_3\\
\end{align*}
while the odd ones are
\begin{align*}
D(\psi^{1,0,n+1}_3)&=\ph^{1,0,m+n+2}_1(m+2)-\ph^{1,1,m+n+1}_3(n+2)
    +\ph^{0,0,n+m+3}_3\\
D(\psi^{0,1,n+1}_3)&=\ph^{0,1,m+n+2}_1(m+2)\\
D(\psi^{0,0,n+2}_1)&=\ph^{0,1,m+n+2}_1(n-m)\\
D(\psi^{0,0,n+2}_2)&=\ph^{1,0,m+n+2}_1(m+2)+\ph^{0,1,m+n+2}_2(n+2)
    +\ph^{0,0,m+n+3}_3\\
D(\psi^{1,1,n}_1)&=\ph^{0,1,m+n+2}_1\\
D(\psi^{1,1,n}_2)&=\ph^{0,1,n+m+2}_2+\ph^{1,1,m+n+1}_3.
\end{align*}

A table of the cohomology is given by
\begin{align*}
(H^1)_o&=\langle \psi^{0,1,0}_3(m+1)+\psi^{0,0,1}_1(m+2)\rangle\\
(H^1)_e&=\langle
\ph^{0,1,0}_1,\ph^{0,1,0}_2(m+2)-\ph^{0,1,0}_2m+\ph^{0,0,0}_3
\rangle\\
(H^{n+2})_o&=\langle \psi^{0,1,n+1}_3-\psi^{1,1,n}_1(m+2),\psi^{0,1,n+1}_3
    +\psi^{0,0,n+2}_1-\psi^{1,1,n}_1(n+2),\\
&\quad-\psi^{1,0,n+1}_3+\psi^{0,0,n+2}_2-\psi^{1,1,n}_2(n+2)
\rangle\\
(H^{n+2})_e&=\langle\ph^{0,1,n+1}_1,\ph^{0,1,n+1}_2+\ph^{1,1,n}_3,\\
&\quad\ph^{0,1,n+1}_2(m+2)+\ph^{0,1,n+1}_2(n+1-m)+\ph^{0,0,n+2}_3
\rangle
\\&\quad\text{if $0\le n<m$}\\
H^{n+2}&=0,\quad \text{if $n\ge m$}.
\end{align*}

Let us analyze the cohomology for $0\le n<m$ first. The first odd
cohomology class in the table is of type $d_{-(n+2)}$ and the last
two are of  type $d_{\#}$.

Since $h_{n+2}=0$, if $n\ge m$, there are no deformations of
$d_{\#}$ on the same level. Thus $d_{\#}$ is rigid in the moduli
space of degree $m+2$ codifferentials.

However, note that we have been so careful up to now to give a
basis for $H_n$ consisting of cohomology classes which correspond
to actual codifferentials, for $n<m$.  But this time, the
cohomology class representing the odd part of $H^1$ doesn't come
from any codifferential of degree 1, and there is no way to fix
this problem, because there is only one cohomology class in degree
1. If one computes the versal deformation of $d_{\#}$, then one
can show that the parameter $t$ which the cocycle by will have to
satisfy $t^2=0$. This means that in some sense, we do not get a
true deformation of $d_{\#}$ in this direction, because we do not
obtain any nonzero values of the parameter for which the
deformation determines an actual codifferential.

\section{Extensions of Codifferentials}

So far, we have not discussed in any detail how to extend a
codifferential $d_N$ which is of fixed degree $N$. Let $G$ be the
subgroup of automorphisms of $S(W)$ fixing $d_N$, and $G_1$ be the
subgroup $G$ consisting of linear automorphisms. The groups $G$
and $G_1$ act on the set of nonzero cohomology classes. We shall
say that two cohomology classes $\delta_1$ and $\delta_2$ are
equivalent if there is an element $f$ in $G$ such that
$f_*(\delta_1)=\delta_2$,  and linearly equivalent if $f$ lies in
$G_1$.

Suppose that $d=d_N + d_{N+k}+ d_{N+k+1}+\cdots$ is an extension
of $d_N$.  The following theorem will help to classify such
extensions.
\begin{thm}\label{equivclasses}
Suppose that $d=d_N + d_{N+k}+ d_{N+k+1}+\cdots$ is a
codifferential. Then $d_{N+k}$ is a cocycle with respect to the
coboundary operator $D=\br{\bullet}{d_N}$. Moreover,
\begin{enumerate}
\item $d$ is equivalent to a codifferential whose first nonzero term
after $d_N$ is of higher order than $N+k$ iff $d_{N+k}$ is a
coboundary,.
\item If $d'=d_N + d'_{N+k}+ d'_{N+k+1}+\cdots$ is another
codifferential and the cohomology class of $d'_{N+k}$ is not linearly
equivalent to the cohomology class of $d_{N+k}$, then $d'$ and $d$
are  not equivalent.
\end{enumerate}
\end{thm}

Let us say that a codifferential $d$ is standard if it is of the
form $d_N + d_{N+k}+ d_{N+k+1}+\cdots$, where $d_{N+k}$ is a
nontrivial cocycle for $d_N$.  By the first part of this theorem,
every nontrivial extension of $d_N$ is equivalent to a standard
codifferential.  For a codifferential in standard form, let us
refer to the cohomology class of $d_{N+k}$ as the secondary term
of $d$.

In general, we don't expect $d_N +d_{N+k}$ to be a codifferential
for an extension $d=d_N + d_{N+k}+ d_{N+k+1}+\cdots$ of $d_N$.
However, for the examples $d_\infty$ and $d_\lambda$ which arise
in this paper, it turns out to be true. We state a theorem which
is useful in characterizing the extensions of a codifferential of
the form $d=d_N +d_{N+k}$.
\begin{thm}\label{equivext}
Suppose $d=d_N +d_{N+k}$ is a codifferential. Let $D_1$ be the
cohomology given by $d_N$, and $D_2$ be the cohomology given by
$d_{N+k}$. Let $d'=d_N +d_{N+k}+d_{N+l}+\cdots$ be an extension of
$d$, where $l>k$. Let $\delta=d_{N+l}$ Then
\begin{enumerate}
\item
$\delta$ is a $D_1$-cocycle. Moreover $D_2(\delta)$ is a
$D_1$-coboundary.
\item
If $\delta$ is a $D_1$-coboundary, then $d'$ is equivalent to an
extension  whose third term
 has degree larger
than $N+l$.
\item If $H^n(D_1)=0$ for $n>N+k$, then any extension $d'$ of $d$ is
equivalent to $d$.
\item
If $\delta=D_2(\eta)$ for some $D_1$-cocycle $\eta$, then $d'$ is
equivalent to an extension whose third term has degree larger than
$N+l$.
\item If every cocycle $\delta$ whose order is larger than $k$
for which $D_2(\delta)$ is a $D_1$-coboundary is of the form
$\delta=D_2(\eta)$ for some $D_1$-cocycle $\eta$, then any
extension of $d$ is equivalent to $d$.
\end{enumerate}
\end{thm}
The proofs of these theorems are in \cite{bfp2}.

\subsection{Extensions of $\d_\lambda$, for $\lambda=-(n+2)$, $n>m$}
The only nontrivial cocycles of degree larger than $m+2$ are of
the form $\psi^{0,0,n+2}_1a$, and for $a\ne 0$, this cocycle is
linearly equivalent to $\psi^{0,0,n+2}_1$. Thus any nontrivial
extension of $d_\lambda$ is of the form
\begin{equation*}
d=\psi^{0,1,m+1}_3-\psi^{1,1,m}_1(n+2)+\psi^{0,0,n+2}_1+\delta+\cdots,
\end{equation*}
where $\delta$ is a cocycle of degree greater than $n+2$, so in
fact, $\delta$ is a coboundary. But then applying
\thmref{equivext}, we see that this term can be eliminated.  Thus
we see that $d$ is equivalent to the extension
\begin{equation*}
\dt{\lambda,e}=\psi^{0,1,m+1}_3-\psi^{1,1,m}_1(n+2)+\psi^{0,0,n+2}_1.
\end{equation*}
Thus there is only one nontrivial extension of  $\dt\lambda$, up
to equivalence.  Let us calculate its cohomology.  Because there
are terms of mixed degree in $\dt{\lambda,e}$,  any cocycle will
also have mixed terms. Let $D_\lambda$, $D_*$, $D_e$ be the
coboundary operators determined by $\dt\lambda$,
$d_*=\psi^{0,0,n+2}_1$, and $\dt{\lambda,e}$, resp. If
$\delta=\delta_{k+2}+\cdots$ is a $D_e$-cocycle of order $k+2$
(order means the least nonzero degree term), then $\delta_{k+2}$
is a $D_\lambda$ cocycle. Moreover,
$D_*(\delta_{k+2})=-D_\lambda(\delta_{n-m+k+2})$. We claim that if
$k>-1$, then the second equality is a consequence of the first
equality.  In fact, we claim that there is an $(n-m+k+2)$-cochain
$\eta$ such that $\gamma=\delta_{k+2}+\eta$ is a $D_e$ cocycle.
Thus $\delta-\gamma$ is a cocycle of higher order than $\delta$.

To see this,  we compute $D_*$ on a basis of $Z_{k+2}(D_\lambda)$,
the space of $D_\lambda$  $(k+2)$-cocycles. Let us suppose that
$k$ is not $-1$ or $m$.
\begin{align*}
&D_*(\psi^{0,1,k+1}_3-\psi^{1,1,k}_1(n+2))=0\\
&D_*(\psi^{1,1,k}_1)=\ph^{0,1,n+k+2}_1
    =D_\lambda(\psi^{0,0,n-m+k+2}_1\tfrac{1}{k-m})\\
&D_*(\psi^{1,0,k+1}_3+\psi^{1,1,k}_2(n-m+k+2))=\\
&\qquad\qquad\ph^{1,0,n+k+2}_1(n+2)+\ph^{0,0,n+k+3}_3
    +\ph^{0,1,n+k+2}_2(n-m+k+2)\\
&\qquad\qquad=D_\lambda(\psi^{0,0,l+2}_2)\\
&D_*(\ph^{0,1,k+1}_1)=0\\
&D_*(\ph^{1,1,k}_3)=-\psi^{1,1,n+k+1}_1(n+2)+\psi^{0,1,n+k+2}_3
    =D_\lambda(-\ph^{0,1,n-m+k+1}_2)\\
&D_*(\ph^{1,0,k+1}_1(n+2)+\ph^{1,0,k+1}_2(k+1-m)+\ph^{0,0,k+2}_3)=0\\
\end{align*}
Now, it is immediate to compute that $\psi^{0,0,n-m+k+2}_1$,
$\psi^{0,0,n-m+k+2}_2$, and $\ph^{0,1,n-m+k+1}$ are all cocycles
with respect to $D_*$.  It follows that every basis element of
$Z_{k+2}(D_\lambda)$ extends to a $D_e$ cocycle.

While it is certainly true that $Z(D_e)$ cannot be decomposed as a
direct product of spaces $Z^n$ which are given by elements of
degree $n$, nevertheless, $Z(D_e)$ is a filtered space and we can
identify $Z^n/Z^{n+1}$ with a subspace of $Z$ in a natural manner.
If $k$ is not $-1$, $m$, or $n$ then $z_{k+2}(D_\lambda)=3|3$ and
the elements below give a corresponding basis for $Z_{k+2}(D_e)$.

\begin{align*}
\eta_1(k)&=\psi^{0,1,k+1}_3-\psi^{1,1,k}_1(n+2))\\
\eta_2(k)&=\psi^{1,0,k+1}_3+\psi^{1,1,k}_2(n-m+k+2)-\psi^{0,0,n-m+k+2}_2\\
\eta_3(k)&=\psi^{1,1,k}_1-\psi^{0,0,n-m+k+2}_1\tfrac{1}{k-m}\\
\tau_1(k)&=\ph^{0,1,k+1}_1\\
\tau_2(k)&=\ph^{1,1,k}_3+\ph^{0,1,n-m+k+1}_2\\
\tau_3(k)&=\ph^{1,0,k+1}_1(n+2)+\ph^{1,0,k+1}_2(k+1-m)+\ph^{0,0,k+2}_3
\end{align*}
One can compute that the $(m+k+3)$-cocycles satisfy
\begin{align*}
\eta_1(m+k+1)&=-D_e(\ph^{0,1,k+1}_2)\\
\eta_2(m+k+1)&=-D_e(\ph^{1,0,k+1,}_2)\\
\eta_3(m+k+1)&=\tfrac{-1}{k+1}D_e(\ph^{1,0,k+1}_1)\\
\tau_1(m+k+1)&=D_e(\ph^{0,1,k+1}_3)\\
\tau_2(m+k+1)&=\tfrac1{n-m+k+2}D_e(-\psi^{1,0,k+1}_3+\psi^{0,0,n-m+k+2}_2)\\
\tau_3(m+k+1)&=D_e(\psi^{0,0,k+2}_2)
\end{align*}
Thus $H^{k+2}=0$, if $k>m$, except possibly for $k=n$.  The reason
we have to single this case out is because then there is an extra
$D_\lambda$ $(n+2)$-cocycle $\eta_4(n)=\psi^{0,0,n+2}_1$, which
remains a cocycle for $D_e$. But $\eta_4(n)=D_e(\ph^{1,0,0}_1)$,
so it is a $D_e$-coboundary, so it does not determine a cohomology
class for $D_e$.

For $k=-1$, one can check that one of the basis elements of
$Z_1(D_\lambda)$ does not extend to a cocycle for $D_e$. We have
\begin{equation*}
Z_1(D_e)=\langle
\ph^{0,1,0}_1,\ph^{1,0,0}(n+2)-\ph^{0,1,0}_2m+\ph^{0,0,1}_3,
\psi^{0,1,0}_3-\psi^{0,0,n-m+1}_1\tfrac{n+2}{m+1}\rangle.
\end{equation*}

Note that when $k=m$,  $\eta_3$ is not defined. The
$D_\lambda$-cocycle $\ph^{1,1,m}_1$, which generates the
cohomology for $D_\lambda$ in dimension $m+2$ does not extend to a
cocycle.  Thus the extension has killed the deformation along the
family. The formula expressing the cocycles of degree $m+k+3$
applies for $k=-1$, except for the $\eta_3$ term, which we don't
need, so all $(m+2)$-cocycles are coboundaries as well.

Thus, we finally arrive at the table of cohomology for
$\dt{\lambda,e}$.
\begin{align*}
H^1&=\langle \psi^{0,1,0}_3(m+1)+\psi^{0,0,n-m+1}_1(n+2),
\\&\quad\ph^{0,1,0}_1,-\ph^{1,0,0}_1(n+2)+\ph^{0,1,0}_2m-\ph^{0,0,1}_3\rangle\\
H^k&=\langle \eta_1,\eta_2,\eta_3,\tau_1,\tau_2,\tau_3\rangle\qquad
    \text{if $0\le k<m$}\\
H^k&=0\qquad\text{if $k\ge m$}
\end{align*}
Notice that the dimension of the cohomology of the extended
codifferential has dropped by one in four different degrees.
This makes sense,  because we have killed nontrivial cocycles in
dimensions 1 and $m+2$,  which create new coboundaries in
dimensions $n+2$ and $m+n+2$, thus killing off two more cohomology
classes.

It is always true that the dimension of the cohomology of an
extended codifferential cannot exceed the dimension of the
cohomology of the original codifferential.  This is because, as in
our example,  cocycles may not extend to cocycles,  but
coboundaries always extend,  since the extended coboundary is
simply the extended codifferential applied to the same object that
gave the original coboundary.

\subsection{Extensions of $\d_\infty$}
The second part of \thmref{equivclasses} guarantees that there are
an infinite number of nonequivalent extensions of $d_\infty=\psi^{1,1,m}_1$.  In
fact, the codifferentials
$\dinfn=\psi^{1,1,m}_1+\psi^{0,1,n+1}_3$, for $n>m$ give an
infinite family of nonequivalent codifferentials extending
$\dt\infty$. Moreover, $\psi^{0,1,n+1}_3a$ is linearly equivalent
to $\psi^{0,1,n+1}_3$ when $a\ne0$, if we can solve the equation
$q^{m-n}=a$.  In particular,  when $\k$ is algebraically closed,
any extension of the form $d=\psi^{1,1,m}_1+\psi^{0,1,n+1}_3a$ is
equivalent to the extension $\dt{\infty,n}$. Let us assume that
$\k$ is algebraically closed.

Let $D_\infty$ be the coboundary operator given by
$\psi^{1,1,m}_1$ and $D_0$ be the one given by $\psi^{0,1,n+1}_3$.
We want to determine how $\dt{\infty,n}$ can be extended.  Since
the next term we add must be a $D_\infty$-cocycle, which may be
taken to be a $D_\lambda$-cohomology class, it may be taken to be
of the form $\psi^{0,1,k+1}_3a$, where $k>n$.  Let
$\zeta(l)=\ph^{0,1,l}_2m-\ph^{0,0,l+1}_3$. We know that $\zeta(l)$
is a $D_\lambda$-cocycle, and moreover
$D_0(\zeta(k-n))=\psi^{0,1,k+l}_3(-m-k+2*n)$, which means that we
can eliminate this term, by \thmref{equivclasses}, unless
$k=2n-m$. However, for this special value of $k$, we cannot
eliminate the term by adding a $D_0$-coboundary of a
$D_\infty$-cocycle, because there is no $D_\infty$-cocycle whose
$D_0$ coboundary coincides with $\psi^{0,1,2n-m+1}_3$. One can
examine the linear equivalences which fix $\dinfn$, and we can
check  that a linear equivalence which fixes $\dinfn$ also fixes
$\psi^{0,1,2n-m+1}_3a$. Thus we cannot eliminate the term by a
linear equivalence, by adding a $D_\infty$-coboundary, or by
adding a $D_0$-coboundary.  This means we cannot eliminate the
term.

Thus we obtain a nontrivial extension of $\dinfn$ given by
\begin{equation*}
\dinfne=\psi^{1,1,m}_1+\psi^{0,1,n+1}_3+\psi^{0,1,2n-m+1}_3a.
\end{equation*}
any extension of $\dt{\infty,n}$ The fact that this is as far as
we need to go is easy to explain.  First,  if we add another term,
it can be taken to be of the form $\psi^{0,1,k+1}_3b$, for some
$k>2n-m$.  But then the term can be eliminated,  because it is a
$D_0$-coboundary.  Notice that the presence of the extra term
$\psi^{0,1,2n-m+1}_3a$ does not interfere with the elimination of
the term by the same process as before.

In fact, it is easy to compute the cocycles for the coboundary
operator $D_e$ induced by $\dt{\infty,n,a}$. For $Z_1$, we note
that the $D_\infty$-cocycles $\ph^{100}_1$, $\ph^{0,1,0}_1$ and
$\psi^{0,1,0}_3$ remain cocycles for $D_a$,  while the cocycle
$\ph^{0,1,0}_2m-\ph^{0,0,1}_3$ does not extend.  This is because
if we let $D_0=\br{\bullet}{\psi^{0,1,n+1}_3+\psi^{0,1,k+1}_3a}$,
we note that
\begin{equation*}
D_0(\ph^{0,1,0}_2m-\ph^{0,0,1}_3)=\psi^{0,1,n+1}_3(n-m).
\end{equation*}
In order for the original cocycle to extend, this must be
$D_\infty$ of something, which it isn't.  It turns out that it
will be more convenient to take a different basis of the cocycles,
which will be more compatible with the basis we will select for
arbitrary degrees.

We obtain the following extensions of $D_\infty$ $(k+2)$-cocycles.
\begin{align*}
\eta_1(k)&=\psi^{0,1,k+1}_3+\psi^{0,1,n-m+k+1}_3\tfrac{a(3n-2m-k)}{2n-m-k},
    \quad k\ne 2n-m\\
\eta_2(k)&=\psi^{1,1,k}_1+\psi^{0,1,n-m+k+1}_3+\psi^{0,1,2(n-m)+k+1}_3a\\
\eta_3(k)&=\psi^{1,1,k}_2-\psi^{1,0,n-m+k+1}_3+\psi^{11,n-m+k}_2(m-k)
\\&-\psi^{1,0,2(n-m)+k+1}_3+\psi^{11,2(n-m)+k}_2a(m-k)\\
\tau_1(k)&=\ph^{1,0,k+1}_1+\ph^{0,1,n-m+k+1}_2(-1+m-k)+\ph^{0,0,n-m+k+2}_3
\\&+\ph^{0,1,2(n-m)+k+1}_2a(-1+m-k)+\ph^{0,0,2(n-m)+k+2}_3a\\
\tau_2(k)&=\ph^{0,1,k+1}_1+\psi^{0,1,n-m+k+1}_1(k-m+1)
\\&+\psi^{0,1,2(n-m)+k+1}_1a(k-m+1)\\
\tau_3(k)&=\ph^{1,1,k}_3+\ph^{1,1,n-m+k}_3(n+1-k+m)
\\&+\ph^{1,1,2(n-m)+k}_3a(n+1-k+m)\\
\end{align*}
When $k=2n-m$ we can use $\psi^{0,1,2n-m+1}_3$ as a basis element
for the $2n-m+2$-cocycles.  What is interesting is that in this
case, it is not a coboundary. The fourth $D_\infty$ even cocycle,
$\zeta(l)$ does not extend to a $D_e$ cocycle, except when
$l=n-m-1$. To see this note that
\begin{equation*}
D_e(\zeta(l))=\psi^{0,1,n+l+2}_3((n-m-l-1)+\psi^{0,1,2n-m+l+2}_3a(2(n-m)-l-1).
\end{equation*}
In general,  we can add  terms of type $\zeta$ to this term in
order to push the degree of the second term up as high as desired,
but the presence
 of the first term prevents us from eliminating
the term altogether. But when $l=n-m-1$, the first term is zero,
so we can add terms to $\zeta(n-m-1)$ in order to create a cocycle
$\tau_4$ of degree $(n-m-1)+2$.  Note that $0\le n-m-1\le n-1$, so
we know something about where that extra even cocycle occurs. The
reason for choosing such a strange basis for the cocycles of
degree $k+2$ is that we can immediately express them as
coboundaries, when the degree is large enough.  In fact,
\begin{align*}
\eta_1(k)&=D_e(\zeta(k-n-1)),\qquad \text{if $k\ge n$, $k\ne 2n-m$}\\
\eta_2(k)&=D_e(-\ph^{0,1,k-m}_2),\qquad\text{ if $k\ge m$}\\
\eta_3(k)&=D_e(\ph^{1,0,k-m}_2),\qquad\text{ if $k\ge m$}\\
\tau_1(k)&=D_e(-\psi^{0,0,k-m+1}_2),\qquad\text{ if $k\ge m$}\\
\tau_2(k)&=D_e(\psi^{0,0,k-m+1}_1),\qquad\text{ if $k\ge m$}\\
\tau_3(k)&=D_e(-\psi^{1,0,k-m}_3),\qquad\text{ if $k\ge m$}\\
\end{align*}
Notice that it is just the corresponding place where
$\zeta(k-n-1)$ extends to a cocycle that we have  an extra odd
cohomology class,  $\eta(2n-m)$.  In fact,  it makes sense that we
should be able to deform in this direction, as we have a family of
extended codifferentials, and this cocycle points in the direction
of the family. Notice that we can deform along the original
family, unlike the previous case where the addition of the
extension term killed the deformation in the direction of the
family.  Of course, in the previous case,  we had that the
extended term was not a cocycle for a deformed element in the
family, but here the extended term is a cocycle for members of the
family, but a trivial one except for $\dt\infty$. Thus as soon as
we deform in the original family direction, the higher order terms
which we added become irrelevant.

Finally,  note that if $m<k<n$,  then we get a deformation in that
direction,  which, according to the moral that things deform to
elements that are earlier, means that we can deform an extension
given by adding a term of degree $n+2$ to an extension given by
adding a term of smaller degree. The cohomology for the extended
codifferential,  as long as $a\ne 0$, is given by
\begin{align*}
H^1&=\langle \eta_1,\tau_1,\tau_2\rangle\\
H^k&=\langle \eta_1,\eta_2,\eta_3,\tau_1,\tau_2,\tau_3\rangle\qquad
    \text{if $0\le k<m$, $k\ne n-m-1$}\\
H^k&=\langle \eta_1 \rangle\qquad \text{if $m\le k<n$}\\
H^k&=0\qquad \text{if $n\le k<2n-m$}\\
H^k&=\langle\eta_1\rangle \qquad\text{if $k=2n-m$}\\
H^k&=0\qquad\text{if $k>2n-m$}
\end{align*}
This completes the description of the extensions of the
codifferentials.
\section{Conclusions}
We have succeeded in giving a complete characterization of all non\-equivalent \linf\ structures
on a $1|2$-dimensional space. The structure of the moduli space of codifferentials starting with a term of degree $N$
is independent of $N$ as long as $N>1$, and moreover, any codifferential of order $N$ generically will deform
to a codifferential of order less than $N$. The deformations of the codifferentials of degree $N$
into codifferentials of the same degree give us a picture of the moduli space of degree $N$ codifferentials.
In two special cases,  degree $N$ codifferentials extend in a nontrivial manner to codifferentials
with terms of higher order, giving us some simple examples of \linf\ structures which are not
equivalent to codifferentials of a single fixed degree.

For the cases which we encountered in this paper,  understanding
the cohomology of the leading order term is sufficient to carry out the complete construction of the
moduli space of codifferentials, and to understand all the extensions.  This situation is actually very
special, compared to the $2|1$-dimensional case, which we studied in \cite{bfp1}.

%%\bibliographystyle{amsplain}
%%\bibliography{global}

\providecommand{\bysame}{\leavevmode\hbox to3em{\hrulefill}\thinspace}
\providecommand{\MR}{\relax\ifhmode\unskip\space\fi MR }
\providecommand{\MRhref}[2]{%
  \href{http://www.ams.org/mathscinet-getitem?mr=#1}{#2}
}
\providecommand{\href}[2]{#2}

\end{document}